\documentclass[twoside,notitlepage,12pt]{article}

\usepackage{amssymb}
\usepackage[leqno]{amsmath}
\usepackage{amsfonts}
\usepackage{amsopn}
\usepackage{amstext}
\usepackage{amsthm}

\renewcommand{\Bbb}{\mathbb}

\newenvironment{pf}{\begin{proof}}{\end{proof}}

\newcommand{\Err}{{\Bbb{R}}}

\newcommand{\al}{\alpha}

\newcommand{\sig}{\sigma}

\renewcommand{\phi}{\varphi}
\renewcommand{\rho}{\varrho}

\newcommand{\rest}{\restriction}

\newcommand{\loe}{\leqslant}
\newcommand{\goe}{\geqslant}

\newcommand{\subs}{\subseteq}

\newcommand{\cl}{\operatorname{cl}}
\newcommand{\Int}{\operatorname{int}}

\newcommand{\id}{\operatorname{id}}

\newcommand{\pr}{\operatorname{pr}}
\newcommand{\liminv}{\varprojlim}
\newcommand{\diag}{\Delta} 

\newcommand{\Es}{{\mathbb{S}}}

\newtheorem{tw}{Theorem}[section]

\newtheorem{lm}[tw]{Lemma}

\theoremstyle{definition}

\theoremstyle{remark}

\newcommand{\setof}[2]{\{#1\colon #2\}}

\newcommand{\pair}[2]{\langle #1, #2 \rangle} 
\newcommand{\map}[3]{#1\colon #2 \to #3} 
\newcommand{\img}[2]{#1[#2]} 
\newcommand{\inv}[2]{{#1}^{-1}[#2]} 

\newcommand{\suppt}{\operatorname{suppt}}

\newcommand{\norm}[1]{\|#1\|}

\newcommand{\invsys}[5]{\langle {#1}_{#4},{#2}_{#4}^{#5},#3 \rangle}

\title{A representation of retracts of cubes}
\author{
{\sc Wies{\l}aw Kubi\'s}\\ \\
Equipe de Logique Math\'ematique\\
Universit\'e Denis-Diderot Paris 7\\
\texttt{wkubis@logique.jussieu.fr}
}

\begin{document}
\maketitle

\begin{abstract}
The purpose of this note is to give the full and self-contained proof of Shche\-pin's result on a spectral representation of retracts of cubes.


\end{abstract}

\section{Preliminaries}

All spaces are assumed to be completely regular. By a {\em cube} we mean any product of compact metric spaces. A {\em map} means, unless otherwise indicated, a continuous map. A {\em retraction} is a map $\map rXY$ which is right-invertible, i.e. there exists a map $\map jYX$ such that $rj=\id_Y$. A particular case is an {\em internal retraction}, i.e. a map $\map rXY$ such that $Y\subs X$ and $r(x)=x$ for every $x\in Y$. In that case one says that $Y$ is a {\em retract} of $X$.

A map $\map fXY$ has a {\em metrizable kernel} if there are a compact metric space $K$ and a map $\map gXK$ such that the diagonal map $f\diag g$ is one-to-one, where $(f\diag g)(x)=\pair{f(x)}{g(x)}$, $x\in X$.

\section{Some lemmas}

Fix a cube $C=\prod_{\al<\kappa}C_\al$ (i.e. each $C_\al$ is a compact metric space) and fix an internal retraction $\map rCX$, i.e. $X\subs C$ and $r(x)=x$ for every $x\in X$. 
A set $S\subs \kappa$ is {\em $r$-admissible} if $x\rest S=x'\rest S$ implies $r(x)\rest S=r(x')\rest S$ for every $x,x'\in C$. Given $S\subs \kappa$, we denote by $\pr_S$ the projection from $C$ onto $C_S=\prod_{\al\in S}C_\al$. Finally, define $X_S=\img{\pr_S}X$ and let $\map{p_S}X{X_S}$ be defined by $p_S(x)=x\rest S$. Assume $S\subs \kappa$ is $r$-admissible. The formula
$$r_S(p_S(x))=r(x)\rest S$$
well defines a map ${r_S}{C_S}{X_S}$, where $C_S=\prod_{\al\in S}C_\al$. Since we deal with quotient maps, $r_S$ is necessarily continuous. If $x\in X_S$ then $x=p_S(y)$ for some $y\in X$ and therefore $r_S(x)=r(y)\rest S=y\rest S=x$. Hence, $r_S$ is a retraction. Later on, we shall show that there are ``many" $r$-admissible sets.

\begin{lm}\label{otwarty} If $S\subs\kappa$ is $r$-admissible then the map $\map{p_S}X{X_S}$ is an open retraction.
\end{lm}

\begin{pf} By the above remarks, $r_S$ is a retraction, which implies that $p_S$ is a retraction, because $p_S r=r_S \pr_S$.

Fix $x_0\in X$ and fix a neighborhood $U$ of $x_0$. It suffices to show that $p_S(x_0)\in\Int \img{p_S}U$. Fix $f\in C(X)$ such that $f\goe0$, $\norm f=1=f(x_0)$ and $\suppt(f)\subs U$, where $\suppt(f)=\cl\setof{x\in X}{f(x)\ne0}$. Define $\map\phi{X_S}{\Err}$ by
$$\phi(x) = \sup \setof{fr(y)}{y\in C,\;\pr_S(y)=x}.$$
In other words, $\phi(x)= \sup \img {fr}{\pr_S^{-1}(x)}$. Since $\pr_S$ is closed and open, $\phi$ is continuous. Clearly $0\loe\phi\loe1$ and $\phi(p_S(x_0))=1$. To finish the proof, we shall show that $\phi(x)=0$ for every $x\in X_S\setminus K$, where $K=\img{p_S}{\suppt(f)}$. Fix $x\in X_S\setminus K$ and fix $y\in C$ such that $\pr_S(y)=x$. Find $g\in C(X_S)$ such that $g\goe0$, $\norm g=1=g(x)$ and $\suppt(g)\cap K=\emptyset$. Let $h=g p_S$. Then $h\in C(X)$ and $\norm{hr+fr}=1$, because $hr$ and $fr$ have disjoint supports. Let $z\in X$ be such that $p_S(z)=x$. Then $z\rest S=x=y\rest S$ and $r(z)=z$. Since $S$ is $r$-admissible, $r(y)=r(z)=z$. Hence $hr(y)=hr(z)=h(z)=g(x)=1$. It follows that $fr(y)=0$, because $fr+hr\goe0$ and $\norm{fr+hr}=1$. This shows that $\phi(x)=0$ and completes the proof.
\end{pf}

The following lemma is trivial.

\begin{lm}\label{trywialny} The union of any family of $r$-admissible sets is $r$-admissible. \end{lm}

The next lemma is crucial.

\begin{lm}\label{manny} Every countable subset of $\kappa$ is contained in a countable $r$-admissible set. \end{lm}

\begin{pf} Fix a countable elementary submodel $M$ of $H(\chi)$, where $\chi>\kappa$ is sufficiently big and $\setof{C_\al}{\al<\kappa},r\in M$. Let $S=M\cap \kappa$. We claim that $S$ is $r$-admissible. For convenience, assume $0\in C_\al$ for each $\al<\kappa$ and therefore there is a sequence $0\in C$ which is constantly equal to $0$. We may assume that this sequence is in $M$. Given $x\in C$ denote by $x\mid S$ the element $x'\in C$ such that $x'(\al)=x(\al)$ for $\al\in S$ and $x'(\al)=0$ for $\al\in\kappa\setminus S$. Now we need to show that $r(x)\rest S=r(x\mid S)\rest S$ for every $x\in C$. Suppose this is not the case and fix $x\in C$ such that $r(x)\ne r(x\mid S)$. Let $u,u_0$ be disjoint basic open sets in $C_S$ such that $r(x)\rest S\in u$ and $r(x\mid S)\rest S\in u_0$. Since each $C_\al$ is second-countable, we may assume that $u,u_0\in M$. Using continuity, find basic open sets $v,v_0\subs C$ such that $x\in v$, $x\mid S\in v_0$ and $\img {\pr_S r}v\subs u$, $\img {\pr_S r}{v_0}\subs u_0$. Then
$$v=\setof{y\in C}{(\forall\;\al\in t)\; y(\al)\in I_\al}$$
and $$v_0=\setof{y\in C}{(\forall\;\al\in r)\; y(\al)\in J_\al},$$
where $t,r$ are finite subsets of $\kappa$ and $I_\al$, $J_\al$ are basic open sets in $C_\al$. Shrinking $v$ and $v_0$ if necessary, we may assume that $t\cap S=r\cap S=s$ and $I_\al=J_\al$ for $\al\in s$. Let $w=\setof{y\in C}{(\forall\;\al\in s)\; y(\al)\in I_\al}$. Then $w\in M$ and $w\cap \inv{(\pr_Sr)}u$ is a nonempty open set (which is witnessed by $x$). Furthermore the set 
$$D=\setof{y\in C}{|\setof{\al<\kappa}{x(\al)\ne0}|\loe\aleph_0}$$ is dense in $C$ and belongs to $M$. By elementarity, there exists $z\in M$ such that $z\in w\cap D$ and $\pr_S(r(z))\in u$. Then the set $\suppt(z)=\setof{\al}{z(\al)\ne0}$ is countable and therefore contained in $M$. Hence $\suppt(z)\subs S$, i.e. $z(\al)=0$ for $\al\notin S$. It follows that $z\in v_0$, which implies $\pr_S r(z)\in u_0$. This is a contradiction, because $u,u_0$ were supposed to be disjoint.
\end{pf}

\section{The representation}

We are now ready to conclude that every retract of a cube is the limit of a ``nice" inverse sequence of open retractions.

\begin{tw} Assume $\setof{C_\al}{\al<\kappa}$ is a collection of metric spaces and $X$ is a retract of $C=\prod_{\al<\kappa}C_\al$. Then $X=\liminv\Es$, where $\Es=\invsys Xp\kappa\al\beta$ is a continuous inverse sequence such that
\begin{enumerate}
	\item $|X_0|=1$;
	\item each $p^{\al+1}_{\al}$ is an open retraction with a metrizable kernel.
\end{enumerate}
\end{tw}

\begin{pf} For each $\al<\kappa$ fix a countable $r$-admissible set $S_\al\subs \kappa$ which contains $\al$ (Lemma \ref{manny}).
Define $A_\al=\bigcup_{\xi<\al}S_\xi$, $X_\al=x_{A_\al}$ and $p_\al=p_{A_\al}$. By Lemma \ref{trywialny}, each $A_\al$ is $r$-admissible. Clearly $X_0$ is a singleton, since $A_0=\emptyset$. For each $\al\loe\beta<\kappa$ there is a unique (necessarily continuous) map $\map{p^\beta_\al}{X_\beta}{X_\al}$ such that $p_\al= p^\beta_\al p_\beta$. By uniqueness we have $p_\al^\beta p_\beta^\gamma=p^\gamma_\al$ for every $\al\loe\beta\loe \gamma$, i.e. $\Es=\invsys Xp\kappa\al\beta$ is an inverse sequence. It is clear that $X=\liminv \Es$ and that the sequence is continuous. Each $p^\beta_\al$ is an open retraction, because, by Lemma \ref{otwarty}, the composition $p^\beta_\al p_\beta=p_\al$ is an open retraction.
Finally, for each $\al<\kappa$, the map $p^{\al+1}_\al$ has a metrizable kernel, because $A_{\al+1}=A_\al\cup S_\al$.
\end{pf}

{\bf Remarks.} One can also get another spectral representation: If $X$ is a retract of a cube, then there exists a $\sig$-complete inverse system of metric spaces $\Es=\invsys Xp\Sigma st$, such that $X=\liminv \Es$ and all the bonding maps as well as all the projections are open retractions. If $X$ is a retract of a Tikhonov cube, then additionally all the bonding maps are soft (see Chapter 2 of Shchepin's article \cite{Szcz}).

\end{document}